\documentstyle{smjrnl} % For Computer Modern Fonts, or,

\newcommand{\be}{\begin{equation}}
\newcommand{\ee}{\end{equation}}

\newcommand{\bd}{\begin{description}}
\newcommand{\ed}{\end{description}}

\begin{document}
\begin{article}

%%%%% To be entered at Kluwers: =====>>
%\journame{}
%\volnumber{}
%\issuenumber{}
%\issuemonth{}
%\volyear{}

%% Do not delete either of the following two commands.
%% Please supply facing curly brackets for the part you
%% are not using for this article.
\received{}\revised{}

\authorrunninghead{Ya.D. SERGEYEV}
\titlerunninghead{Parallel information algorithm with local tuning}

\pagestyle{empty}

%%%%  <<==== End of commands to be entered at Kluwers

%%  Authors, start here ====>>

\title{Parallel Information Algorithm with Local Tuning \
for Solving Multidimensional GO Problems}

%% Author name
\authors{Yaroslav D. Sergeyev}
%% Email address. Do not delete this command even if you do not have
%% an email address; just leave the command followed by facing curly brackets.
\email{yaro@si.deis.unical.it}

%% Affiliation address
\affil{Institute of Systems Analysis and Information Technology of National Research Council, c/o DEIS, University of Calabria, 87036 Rende (CS), Italy and\ Software Department, University of Nizhni Novgorod, Gagarin Av. 23, Nizhni Novgorod, Russia}

%% Article editor
%%\editor{}

%% Abstract
\abstract{In this paper we propose a new parallel algorithm for solving global optimization (GO) multidimensional problems. The method unifies two powerful approaches for accelerating the search : parallel computations and local tuning on the behavior of the objective function. We establish convergence conditions for the algorithm and theoretically show that the usage of local information during the global search permits to accelerate solving the problem significantly. Results of numerical experiments executed with 100 test functions are also reported.}

%% Keywords

\section{Introduction}
Let us consider the problem of finding the global minimum and global minimizers of a function $\phi (y), y \in  D$, where
\[
D =  \{ x : a_{i} \le x_{i} \le b_{i}, 1\le i\le N \}
\]
is a hyperinterval in $R^{N}$. This problem is actively studied by many authors (see e.g. \cite{Floudas and Pardalos (1996)}, \cite{Horst and Pardalos (1995)}, \cite{Migdalas Pardalos and Storoy (1997)}, \cite{Pardalos Phillips and Rosen (1992)}, \cite{Strongin (1978)}). One of the promising approaches to attack this problem is in using the Peano type space-filling curves.

It is known (see  \cite{Butz (1968)},  \cite{Strongin (1978)},  \cite{Strongin (1992)}) that this multidimensional problem can be reduced to a one-dimensional one by using the curves. In this case we obtain
\begin{equation}
 \min\{\phi (y): y\in D\} =  \min\{\phi (y(x)): x\in [0,1]\},   \label{(3.17)}
\end{equation}
moreover, if the multidimensional function $\phi (y)$ satisfies  the Lipschitz condition with the constant $L$ over $D$  then  $\phi (y(x))$ over the interval $[0,1]$ satisfies the H\H{o}lder condition
\begin{equation}
\mid \phi (y(x')) - \phi (y(x''))\mid  \le  H\mid x' - x''\mid ^{1/N}, \hspace{5mm} x', x'' \in  [0,1],      \label{(3.18)}
\end{equation}
where $H$ is the H\H{o}lder constant and
\[
 H = 4Ld\sqrt{N}, \qquad d = \max \{ b_{i}-a_{i}: 1\le i\le N \}.
\]
The constant $H$ corresponds to the  space-filling  curve $y(x)$ which computational scheme is presented  in  \cite{Strongin (1978)} (see also \cite{Strongin and Sergeyev (1992)}). A special procedure for fast calculating the image $y(x) \in D$ on the curve for every given $x \in  [0,1]$ has been proposed in these papers. The procedure has a high speed because, instead of constructing an approximation of the whole curve, it  directly computes the coordinates of $y(x)$ in $D$ for every given $x$. In addition, its parallel version has been introduced in \cite{Strongin and Sergeyev (1992)}.

For solving the problem (\ref{(3.17)}), (\ref{(3.18)})  the  {\it information algorithm} (IA) has been proposed in \cite{Strongin (1978)}, \cite{Strongin (1992)} for the case when the Lipschitz constant (and,  therefore  the  H\H{o}lder  one)  is unknown. A special procedure for estimating $H$ on the base  of information obtained during the  search  has  been elaborated.

The power of parallel computations is widely used in global optimization  to accelerate the search  (different approaches are presented in \cite{Bertocchi (1990)}, \cite{Grishagin Sergeyev and Strongin (1997)},  \cite{Migdalas Pardalos and Storoy (1997)} -- \cite{Schendel (1984)}, \cite{Sergeyev and Grishagin (1994a)}, \cite{Sergeyev and Grishagin (1994b)},  \cite{Strongin and Sergeyev (1992)} -- \cite{Van Laarhoven (1985)}). It has been shown that very often simple parallelizing ideas (parallel grid method or subdividing the search region in $p$ subregions where $p$ processors work in parallel) can lead to appearence of redundant  function evaluations in comparison with an efficient sequential method. That is why one of the  reasonable ways for  usage of parallel computers consists in parallelizing fast sequential methods.

Wide numerical experiments and a deep theoretical study have shown (see \cite{Grishagin (1978)}, \cite{Strongin (1978)}) that IA has good speed characteristics in comparison with other methods which don't use derivatives. In \cite{Strongin and Sergeyev (1992)} a {\it parallel information algorithm} (PIA) extending IA to the case of parallel computations has been introduced. Convergence conditions and estimates of speed up obtained in comparison with the original sequential IA have been established. A special study has been executed to obtain conditions which guarantee absence of redundant (in comparison with the sequential case) evaluations of $\phi (y(x))$. This idea has demonstrated to be fruitful and it was applied to parallelize some other sequential GO methods (see \cite{Gergel and Sergeyev}, \cite{Sergeyev and Grishagin (1994a)}, \cite{Sergeyev and Grishagin (1994b)}). By generalizing this approach a class of parallel characteristical GO algorithms has been introduced and theoretically investigated in \cite{Grishagin Sergeyev and Strongin (1997)}.

An alternative approach for acceleration of the global search consists in the following. It has been shown  for different one-dimensional GO algorithms (see \cite{Sergeyev (1995a)} -- \cite{Sergeyev (1998)}) that using  estimates of the Lipschitz constant $L$ slows down the search for that subregions where the local Lipschitz constants are significantly less than $L$ (hereinafter we shall call $L$ the  {\it global} Lipschitz constant).  The sequential {\it information} GO {\it algorithm with local tuning} (IALT) has been introduced in \cite{Sergeyev (1995b)} for solving the problem (\ref{(3.17)}), (\ref{(3.18)}). In that paper a special procedure has been proposed to estimate adaptively  local  constants in different areas of the search region to accelerate the search. Theoretical results and numerical experiments have shown that this approach permits to accelerate the search significantly in comparison with the original IA using esimates of the global Lipschitz constant.

In this paper we propose to unify both approaches and construct a {\it parallel} information algorithm with {\it local tuning} (PLT). In order to accelerate the search during every iteration the new method makes the following :
\bd
\item
- it adaptively estimates {\it local H\H{o}lder constants} in different subintervals  of $[0,1]$ to tune itself  on the local behaviour of the reduced objective function $\phi (y(x)), x \in  [0,1]$;
\item
- it simultaneously evaluates the objective function $\phi (y(x))$ at $p$ different points in $p$ intervals (having the highest probability to contain the global optimum) on $p$ parallel computers;
\item
- local information about $\phi (y(x))$ is used over the {\it whole} search region {\it during} the global search in contrast with traditional approaches (see e.g. \cite{Floudas and Pardalos (1996)}, \cite{Horst and Pardalos (1995)}) usually starting a local search in a neighbourhood of a global minimizer after stopping the global procedure.
\ed

It will be shown further that  adaptive estimates of the local H\H{o}lder constants introduced here permit to accelerate the search significantly in comparison with the parallel methods  (see \cite{Grishagin Sergeyev and Strongin (1997)}, \cite{Strongin and Sergeyev (1992)}) using the global Lipschitz (H\H{o}lder) constants or their estimates.

The rest of the paper is structured as follows. Section 2 describes the computational scheme of PLT and presents convergence results of the new type for the introduced parallel method. Section 3 contains results of a numerical comparison between IA, PIA, IALT and PLT executed with $100$ two-dimensional functions usually used in literature for numerical comparison between the information GO methods (see \cite{Grishagin (1978)}, \cite{Grishagin Sergeyev and Strongin (1997)}, \cite{Strongin (1978)}, \cite{Strongin and Sergeyev (1992)}). Some conclusions are presented in the last Section.

\section{The algorithm and its convergence conditions}

To describe the algorithm we need some designations and definitions. We introduce a simplifying designation
\[
f(x) = \phi (y(x)), \hspace{5mm} x \in  [0,1],
\]
for the objective function $\phi (y(x))$. We call {\it trial} the operation of evaluating $f(x)$ at a point $x$. The character $l$ we use as the counter of parallel iterations of PLT. In the course of
every parallel iteration $l$ we evaluate $f(x)$ at $p(l)>1$ points using $p(l)$ parallel processors. We use the Peano curve approximation with the depth $m$ presented in \cite{Strongin (1978)}, \cite{Strongin and Sergeyev (1992)} and a preset accuracy $\epsilon $ concorded with $m$ by the inequality
\[
\epsilon  \ge  2^{-m}/(4\sqrt{N}).
\]
Let us now present the computational scheme of the new method.
\begin{description}
\item
{\bf Step 0.} Execute $q(l) > 1 $ initial trials at the points $x^{1}= 0, x^{2}= 1$ and  some  internal  points  $x^{3}, x^{4}, \ldots   ,x^{q(l)}$ belonging to $(0,1)$. Set $l = 1$.
\item
{\bf Step 1.} Points $x^{1},x^{2},\ldots  , x^{q(l)}$ of the previous trials reorder by increasing of their coordinates using subscripts, i.e.
\[
0 = x_{1} < x_{2} < \ldots   < x_{q-1}< x_{q}= 1
\]
where $q = q(l)$.
\item
{\bf Step 2.} Calculate the values $\mu _{j}$ estimating local H\H{o}lder constants for the intervals $[x_{j-1},x_{j}],  2\le j\le q,$ following the rules :
\be
\mu _{j}= \max \{ \lambda _{j}, \gamma _{j}, \xi \},\; 2 \le  j \le  q, \label{loc1}
\ee
where $\xi  > 0$ is a small number, the value $\lambda _{j}$ spying on the local information is determened by the formulae
\be
\lambda _{j}= \max \{ \frac {\mid f(x_{i})-f(x_{i-1}) \mid}
                  {\mid x_{i}-x_{i-1}\mid ^{1/N}} : i \in I_j \},  \label{loc2}
\ee
\[
I_j =      \left\{ \begin{array}{ll}
                \{2,3\}          & \mbox{ if $j=2$}             \\
                \{j-1,j,j+1\}    & \mbox{ if $3 \le j \le q-1$} \\
                \{q-1,q\}        & \mbox{ if $j=q$}
               \end{array}
     \right \}.
\]
The last part of formulae (\ref{loc1}) controls the global information accumulated during all the previous iterations and is represented by the value
\be
\gamma _{j}= \mu (x_{j}- x_{j-1})^{1/N}/(X^{\max})^{1/N}, \label{loc3}
\ee
where
\[
\mu = \max \{\mid f(x_{j})-f(x_{j-1})\mid /(x_{j}-x_{j-1})^{1/N}: 2\le j\le q\}, \label{loc4}
\]
\[
X^{\max}= \max \{ x_{i}- x_{i-1}: 2 \le  i \le q \}.
\]
\item
{\bf Step 3.} For the intervals $[x_{j-1}, x_{j}], 2\le j\le q,$ compute their characteristics $R(j)$ as follows
\be
\hspace{-0.9cm}R(j) = r\mu _{j}(x_{j}-x_{j-1})^{1/N}+\frac {(z_{j}-z_{j-1})^{2}}
                         {r\mu _{j}(x_{j}-x_{j-1})^{1/N}}-(z_{j}+z_{j-1}), \quad 2\le j\le q,  \label{R}
\ee
where $z_{j}=f(x_{j}), 1 \le j \le q,$ and $r > 1$ is the reliability parameter of the method.
\item
{\bf Step 4.} Calculate the new trial points $x^{q+j}\in  (0,1), 1 \le j\le p$, and their images $y(x^{q+j}) \in  D, 1 \le j\le p$, on the Peano curve $m$-approximation
\begin{equation}
\hspace{-0.9cm}x^{q+j} = 0.5(x_{t_{j}-1}+ x_{t_{j}}) -
         (2r)^{-1}(\mid z_{t_{j}}-z_{t_{j}-1}\mid / \mu_{t_{j}})^{N}\mbox{sign}(z_{t_{j}}-z_{t_{j}-1}),            \label{(4)}
\end{equation}
where
\[
t_{1}=\mbox{argmax} \{ R(i): 1 < i \le  q \},
\]
\[
t_{j}=\mbox{argmax} \{ R(i): 1<i\le q, i\neq t_{s}, 1\le s\le j-1\}, 1<j\le p,
\]
\[
p = p(l+1) \le  q(l)-1.
\]
\item
{\bf Step 5.} Execute trials of the $(l+1)$th iteration at the points $y(x^{q+j}), 1 \le j\le p$, where $x^{q+j}, 1 \le j\le p$, are from (\ref{(4)}).
\item
{\bf Step 6.} If the stopping rule
\[
\min \{ \mid x_{t_{j}}-x_{t_{j}-1}\mid ^{1/N}: 1  \le  j \le  p \} \le  \epsilon
\]
is not satisfied, then go to Step 1. Otherwise take the value
\[
z^{*}_{q}  = \min \{ z_{i}: 1 \le  i \le  q \}
\]
as an estimate of the global minimum and corresponding to this value points $y^{*}_{q} \in D$ as estimates of global minimizers and Stop.
\end{description}

The algorithm described above belongs to the class of parallel characteristical methods introduced in \cite{Grishagin Sergeyev and Strongin (1997)}. In comparison with the other algorithms from the class \cite{Grishagin Sergeyev and Strongin (1997)} the main peculiarity of the new method consists in the following. Instead of usage the global H\H{o}lder constant $H$ or its estimates, PLT executes an adaptive local tuning on the behaviour of the reduced one-dimensional function $f(x)$. The local tuning is based  on the  estimates $\mu _{j}$ of the local H\H{o}lder constants $H_j$ for every subinterval $[x_{j-1},x_{j}],  2\le j\le q(l),$ during every iteration $l > 1$.

The value $\mu _{j}$ from (\ref{loc1}) represents the result of a balance between the global and local information accumulated during the search. If the interval $[x_{j-1}, x_{j}]$ is small then, the local information (see (\ref{loc2})) represented by $\lambda _{j}$ is very important. The global information (see (\ref{loc3})) represented by $\gamma _{j}$ is less important because the estimate $\mu$ of the global H\H{o}lder constant (see (\ref{loc4})) could be obtained at an interval being very far from $[x_{j-1}, x_{j}]$. Vice versa, when the interval $[x_{j-1}, x_{j}]$ is   wide, then we cannot trust the local information from (\ref{loc1}) and $\gamma _{j}$ plays the most important part. The parameter $\xi $ reflects our idea  that $f(x)$ is such that over every subinterval $[x_{j-1}, x_{j}]$  its  H\H{o}lder constant $H_j \ge \xi$.

Let us now consider convergence properties of the new method. Let $\{x^{q}\}$ be the sequence of the trial points generated by PLT during minimizing the function $f(x)$. The first theorem on a level with  other results asserts that only local minimizers of $f(x)$ can be limit points  of $\{x^{q}\}$.

\vspace{3mm}
{\bf Theorem 1.} Let $\bar x  \in  [0,1]$ be a limit point of  the sequence $\{x^{q}\}$ and
\begin{equation}
p(l) \le  Q < \infty , \; l>1.   \label{(3.28)}
\end{equation}
Then the following results take place :
\begin{description}
\item
{i.} if $\bar x \in  (0,1)$ then, there exist two subsequences of $\{x^{q}\}$ such that the first one converges to $\bar x $  from the left and the second one - from the right;
\item
{ii.} the point $\bar x $ is a local minimizer of $f(x)$ if the function $f(x)$ has a finite number of local extrema;
\item
{iii.} if a part $\bar x $ there exists another limit point $\hat{x}$, then $f(\bar x ) = f(\hat{x})$;
\item
{iv.} for all $q \ge  1$ it follows $z^{q}= f(x^{q}) \ge  f(\bar x )$.
\end{description}

{\bf Proof.} Since PLT belongs to the class of parallel characteristical methods these results can be easily deduced from the general convergence theory presented in \cite{Grishagin Sergeyev and Strongin (1997)}. \rule{5pt}{5pt}

The second theorem describes sufficient convergence conditions of the sequence $\{x^{q}\}$ to a global minimizer $x^{*}$. This result can not be obtained within the general framework from \cite{Grishagin Sergeyev and Strongin (1997)} and we prove it. In order to proceed let us introduce $\{q\} =  \{1,2,3,\ldots   \}$ as the sequence enumerating iterations executed by PLT.

\vspace{3mm}
{\bf Theorem 2.} If (\ref{(3.28)}) is true and there exists an infinite subsequence $\{h\} ,  \{h\} \subseteq  \{q\}$, such that for the interval $[x_{j-1}, x_{j}], j=j(l), l\in \{h\}, $ containing the point $x^{*}$ during the $l$th PLT iteration the following inequality takes place
\be
r\mu _{j}\ge  2^{1-1/N}K_{j}+ (4^{1-1/N}K^{2}_{j}- M^{2}_{j})^{1/2}, \label{28}
\ee
where
\be
K_{j}= \max \{ (z_{j-1}- f(x^{*}))(x^{*} - x_{j-1})^{-1/N}, (z_{j}- f(x^{*}))(x_{j}- x^{*})^{-1/N}\}, \label{29}
\ee
\be
M_{j}= \mid z_{j-1}- z_{j}\mid (x_{j}- x_{j-1})^{-1/N}, \label{30}
\ee
then, $x^{*}$ is a limit point of $\{x^{q}\}$.

{\bf Proof.} In order to start the proof we notice that the following inequality takes place for the estimates $\mu _{j}$ from (\ref{loc1}):
\be
0 < \xi \leq \mu _{j(l)} \leq  \max \{ H, \xi\} < \infty, \hspace{1cm} 2 \le  j \le  q(l),  l \geq 1.
\label{23}
\ee

Suppose now that there exists a limit point $x' \neq x^{*}$ of the trial sequence $ \{ x^{q} \}$. From (\ref{R}), (\ref{23}) and the first assertion of Theorem 1 we conclude for the interval $[ x_{i-1}, x_{i} ], i =i(l)$, containing $x'$ at the $l$th iteration of PLT, that
\be
    \lim_{l \rightarrow \infty} R(i(l))=-4f(x').  \label{16}
\ee
Consider now the interval $[ x_{j-1}, x_{j}], j=j(l)$,
\be
    x^{*} \in [x_{j-1}, x_{j}]  \label{17}
\ee
and suppose that $x^{*}$ is not a limit point of $\{x^{q} \}$. This means that there exists an iteration number $m$ such that for all $l \geq m$
\[
x^{q(l)+k} \notin[ x_{j-1}, x_{j}], \hspace{1cm} j =j(l),  1 \le k\le p(l+1),
\]
i.e. new trial points will not fall in the interval (\ref{17}).
Estimate now the characteristic $R(j(l)), l \geq m,$ of this interval.
    It follows from (\ref{29}) and (\ref{17}) that
\[
z_{j-1} - f(x^{*}) \leq K_{j} ( x^{*} - x_{j-1} )^{1/N},
\]
\[
z_{j} - f(x^{*}) \leq K_{j} ( x_{j} - x^{*} )^{1/N}.
\]

By summrizing these two inequalities and by using  the designation
\[
\alpha = (x^{*} - x_{j-1} ) / ( x_{j} - x_{j-1} ),
\]

we obtain
\[
z_{j-1} + z_{j} \leq 2f(x^{*}) + K_{j} ( ( x^{*} - x_{j-1} )^{1/N} + ( x_{j} - x^{*} )^{1/N} ) =
\]
\[
= 2f( x^{*} ) + K_{j} ( \alpha^{1/N} + ( 1 - \alpha )^{1/N} ) ( x_{j} - x_{j-1} )^{1/N} \leq
\]
\[
 \leq 2f( x^{*} ) + K_{j} ( x_{j} - x_{j-1} )^{1/N} max \{ \alpha^{1/N} + (1 - \alpha )^{1/N} : 0 \leq a \leq 1 \} =
\]
\[
= 2 ( f( x^{*} ) + 2^{-1/N} K_{j} ( x_{j} - x_{j-1} )^{1/N} ).
\]
From this estimate,  (\ref{28}), and (\ref{30}) we have
\[
R(j(l)) = r \mu_{j} ( x_{j} - x_{j-1} )^{1/N} + ( z_{j-1} - z_{j} )^{2} ( r \mu_{j} )^{-1} (x_{j} - x_{j-1} )^{-1/N} - 2( z_{j-1} + z_{j} )
\]
\be
 \geq ( x_{j} - x_{j-1} )^{1/N} (r \mu_{j} + M^{2}_{j} (r \mu_{j} )^{-1} - 2^{2-1/N} K_{j} ) - 4 f(x^{*} ) \geq - 4 f(x^{*} )   \label{20}
\ee
for all iteration numbers $l \in \{ h \} $.

Since $x^{*}$ is a global minimizer and the sequence $ \{h\}$ is infinite, then from (\ref{16}) and (\ref{20}) it follows that an iteration number $l^{*}$ will exist such that

\be
R(j(l^{*})) \geq R(i(l^{*})).  \label{21}
\ee
But, in according with the decision rules of PLT, this means that during the $l^{*}$th iteration one of $p(l^{*})$ new trials will be executed at the interval (\ref{17}).  Thus, our assumption that $x^{*}$ is not a limit point of $\{x^{q} \}$ is not true and theorem has been proved.
 \rule{5pt}{5pt}

\vspace{3mm}
    {\bf Corollary 1.} Given the conditions of Theorem 2 all the limit points  of the sequence $\{x^{q} \}$ are the global minimizers of $f(x)$.

    {\bf Proof.} The corollary follows from the third assertion of Theorem 1.   \rule{5pt}{5pt}

    Let $X^{*}$ be the set of the global minimizers of the function $f(x)$. Corollary 1 ensures that the set of limit points of the sequence $\{x^{q} \}$ belongs to $X^{*}$. Conditions ensuring coincidence of  these sets are  established by Corollary 2.

\vspace{3mm}
    {\bf Corollary 2.} If condition (\ref{17}) is fulfilled for all the points  $x^{*} \in X^{*}$, then the set of limit points of $\{x^{q} \}$ coincides with $X^{*}$.

    {\bf Proof.} Again, the corollary is a straightforward consequence of Theorem 2 and the third assertion of Theorem 1.   \rule{5pt}{5pt}

Theorem 2 and its corollaries are very important both from theoretical and practical viewpoints. It is known (see   \cite{Butz (1968)},  \cite{Strongin (1978)},  \cite{Strongin (1992)}) that every point $y \in D$ can have up to $2^N$ images  on the curve. The global minimizer $ x^{*} $ can have up to $2^N$ images on the curve too. To obtain an $\epsilon$-approximation of $ x^{*} $ it is enough to find only {\it one} its image on the curve.

We have proved that to have convergence to a global minimizer $x^{*}$ PLT needs  the fulfilment of  condition (\ref{28}) (which is considerably weaker than the Lipschitz condition) for one of the images of the point $ x^{*} $ on the curve. Thus, (in contrast with the other methods from the class \cite{Grishagin Sergeyev and Strongin (1997)}) the new parallel method does not need  the exact value of the precise Lipschitz constant $L$ (neither its upper estimate) for the {\it whole} region $D$. It is enough that condition (\ref{28}) is fulfilled in a neighbourhood of $ x^{*} $ for {\it one} its image on the curve. In contrast with this, methods using in their work the exact  Lipschitz (H\H{o}lder) constant  (or its upper estimate)  will have convergence to {\it all  $2^N$ images} on the curve. Of course, this fact leads to a significant slowing down the search and explains why PLT works faster. Results  of numerecal experiments to be presented in the next Section confirm that the new approach permits to accelerate the search significantly.

\section{Numerical experiments}
In this section we compare performance of the new method PLT with the original information algorithm  (IA), parallel information algorithm (PIA), and sequential information algorithm with local tuning (IALT). An ALLIANT FX/80 parallel mini-supercomputer having 4 processors has been used on the series of 100 two-dimensional multiextremal functions from \cite{Grishagin (1978)} usually applied for testing information GO methods (see \cite{Grishagin (1978)}, \cite{Grishagin Sergeyev and Strongin (1997)}, \cite{Sergeyev (1995b)}, \cite{Strongin (1978)}, \cite{Strongin and Sergeyev (1992)}):
\begin{equation}
\hspace{-8mm} f(x) = \sqrt { \Bigl[ \sum^{7}_{i=1} \sum^{7}_{j=1} (A_{ij}a_{ij}(x)+B_{ij}b_{ij}(x))\Bigr]^{2} +
 \Bigl[ \sum^{7}_{i=1} \sum^{7}_{j=1} (C_{ij}a_{ij}(x)-D_{ij}b_{ij}(x))\Bigr]^{2} },  \label{(4.1)}
\end{equation}
where $0 \le  x_{1} \le  1, 0 \le  x_{2} \le  1$ and
\[
a_{ij}(x) = \sin (i\pi x_{1}) \sin (j\pi x_{2}),\qquad
b_{ij}(x) = \cos (i\pi x_{1}) \cos (j\pi x_{2}),
\]
and $A_{ij},  B_{ij}, C_{ij}, D_{ij}$ are random coefficients from the interval [-1,1].

For all the methods in all the experiments we have used the 12-order approximation of the Peano curve, initial points
$\{ 0.2, 0.4, 0.6, 0.9 \}$, the reliability parameter $r = 2.9$ and the search accuracy $\epsilon  = 0.001$. We have chosen $\xi  = 10^{-6}$ in PLT.

In Tables $1, 2$ we present average results for 100 functions from the class (\ref{(4.1)}). The column "\%" shows a quantity of experiments in which global minima have been found. In Tab. 1 we compare the sequential method IALT and PLT. It is seen from the table that the introduced type of the parallelizm permits to achive high levels of speed up in comparison with the sequential method IALT. Note, that obtaining speed up higher than the number of processors used in the parallel case is possible due to parallel adaptive estimating the local H\H{o}lder constants.

In order to underline the effect obtained after introducing the local tuning, in Tab. 2 we report the values of speed up obtained by using the parallel method PIA in comparison with the original sequential IA (see \cite{Strongin and Sergeyev (1992)}). Both methods use an adaptive estimate of the  {\it global} H\H{o}lder constant in their work.

In Tab.$3$ we compare PLT working with the local H\H{o}lder constants with PIA using the global H\H{o}lder constant. The data from Tabs. $1, 2$ have been used. The obtained values of speed up both in time and in trials are shown. It is seen from Tab. $3$ that in comparison with PIA the new method PLT functions faster more than $4.5$ times in trials and more than $3.5$ times in time.

\begin{table}[ht]
\caption{Average results of the numerical experiments executed by PLT w.r.t. the sequential method IALT.}
\begin{tabular*}{\textwidth}{@{\extracolsep{\fill}}ccccccc}
\hline
Method&Processors&\%&Trials&Time&Speed up&Speed up\cr
 & & & & & (trials)& (time)\cr
\hline
IALT &   1&98&351.37&11.15&-&-\cr
PLT&   2&96&339.00&4.22&1.97&3.33\cr
PLT&   3&98&349.75&2.73&3.02&7.15\cr
PLT&   4&98&348.24&1.99&3.94&10.33\cr
\hline
\end{tabular*}
\end{table}

\begin{table}[ht]
\caption{Average results of numerical experiments with the information algorithms using global estimates of the H\H{o}lder constant.}
\begin{tabular*}{\textwidth}{@{\extracolsep{\fill}}ccccccc}
\hline
Method&Processors&\%&Trials&Time&Speed up&Speed up\cr
 & & & & & (trials)& (time)\cr
\hline
IA &   1&100&1575.12&70.04&-&-\cr
PIA&   2&100&1596.08&21.06&1.97&3.33\cr
PIA&   3&100&1562.61&9.80&3.02&7.15\cr
PIA&   4&100&1599.92&6.78&3.94&10.33\cr
\hline
\end{tabular*}
\end{table}

\begin{table}[ht]
\caption{Speed up obtained by using PLT in comparison with PIA on the data from Tabs. $1, 2$.}
\begin{tabular*}{\textwidth}{@{\extracolsep{\fill}}ccc}
\hline
Processors&Speed up in trials&Speed up in time\cr
\hline
2&4.71&4.99\cr
3&4.47&3.59\cr
4&4.59&3.41\cr
\hline
\end{tabular*}
\end{table}

\section{A brief conclusion}

In this paper a new parallel algorithm for solving global optimization  multidimensional problems has been proposed. The method unifies two powerful approaches for accelerating the search : parallel computations and local tuning on the behavior of the objective function. Peano-type space-filling curves have been used to reduce the multidimensional problem to the one-dimensional one. For the obtained problem the new method  adaptively estimates local H\H{o}lder constants over different subintervals  of the one-dimensional search region to tune itself  on the local behaviour of the reduced objective function. The local information  is used by the method over the {\it whole} search region {\it during} the global search in contrast with traditional approaches  usually starting a local search in a {\it neighbourhood} of a global minimizer {\it after} stopping the global procedure. It has been theoretically and  numerically shown that the new technique permits to accelerate the search considerably.

\acknowledgements
{This research was partially supported by the Russian Foundation of Fundamental Research, grant 95-01-01073.}

\end{article}

\begin{references}

\bibitem{Bertocchi (1990)}Bertocchi, M. (1990), A parallel algorithm for global optimization, {\it Optimization}, {\bf 62(3)}, 379--386.

\bibitem{Butz (1968)}Butz, A.R. (1968), Space filling curves and mathematical programming, {\it Inform. Control.}, {\bf 12(4)}, 314--330.

\bibitem{Gergel and Sergeyev}Gergel, V.P. and Sergeyev, Ya.D. Sequential and parallel global optimization algorithms using derivatives, (to appear in {\it
          Computers and Mathematics with Applications}).

\bibitem{Grishagin (1978)}Grishagin, V.A. (1978), Operation characteristics of some global optimization algorithms, {\it Problems of Stochastic Search}, {\bf 7}, 198--206.

\bibitem{Grishagin Sergeyev and Strongin (1997)}Grishagin, V.A., Sergeyev, Ya.D. and Strongin, R.G. (1997),  Parallel characteristical global optimization algorithms, {\it J. of Global Optimization}, {\bf 10}, 185--206.

\bibitem{Floudas and Pardalos (1996)}Floudas, C.A. and Pardalos, P.M. (1996), Eds., {\it State of the Art in Global Optimization}, Kluwer Academic Publishers, Dordrecht.

\bibitem{Horst and Pardalos (1995)}Horst, R. and Pardalos, P.M. (1995), Eds., {\it Handbook of Global Optimization}, Kluwer Academic Publishers, Dordrecht.

\bibitem{Migdalas Pardalos and Storoy (1997)}Migdalas, A., Pardalos, P.M., and Stor\o{}y, S. (1997), Eds., {\it Parallel Computing in Optimization}, Kluwer Academic Publishers, Dordrecht.

\bibitem{Pardalos Phillips and Rosen (1992)}Pardalos, P.M., Phillips, A.T., and Rosen, J.B. (1992), {\it Topics in Parallel Computing in Mathematical Programming}, Science Press.

\bibitem{Schendel (1984)}Schendel, U. (1984), {\it Introduction to Numerical Methods for Parallel Computers}, Ellis Horwood, Chichester, England.

\bibitem{Sergeyev (1995a)}Sergeyev Ya.D. (1995) A one-dimensional deterministic global minimization algorithm,  {\it Comput. Maths. Math. Phys.}, {\bf 35(5)},
     705--717.

\bibitem{Sergeyev (1995b)}Sergeyev, Ya.D. (1995), An information global optimization algorithm with local tuning,  {\it SIAM J. Optimization}, {\bf 5(4)}, 858--870.

\bibitem{Sergeyev (1998)}Sergeyev, Ya.D. (1998), Global one-dimensional optimization using smooth auxiliary functions, {\it Mathematical Programming},
   {\bf  81(1)}, 127--146.

\bibitem{Sergeyev and Grishagin (1994a)}Sergeyev, Ya.D., and Grishagin, V.A. (1994), A parallel method for finding the global minimum of univariate functions, {\it JOTA}, {\bf 80}, 513--536.

\bibitem{Sergeyev and Grishagin (1994b)}Sergeyev, Ya.D., and Grishagin, V.A. (1994), Sequential and parallel algorithms for global optimization, {\it Optimization Methods and Software}, vol.3, 111-124.

\bibitem{Strongin (1978)}Strongin, R.G. (1978), {\it Numerical Methods on Multiextremal Problems}, Nauka, Moscow.

\bibitem{Strongin (1992)}Strongin, R.G. (1992),  Algorithms for multiextremal mathematical programming problems employing the set of joint space-filling curves, {\it J. of Global Optimization}, {\bf 2}, 357--378.


\bibitem{Strongin and Sergeyev (1992)}Strongin, R.G. and  Sergeyev, Ya.D. (1992), Global multidimensional optimization on parallel computer, {\it Parallel Computing}, {\bf 18},  1259--1273.

\bibitem{Sutti (1984)}Sutti, C. (1984), Local and global optimization by parallel algorithms for MIMD systems, {\it Annals of Operations Research}, {\bf 1}, 151--164.

\bibitem{Van Laarhoven (1985)}Van Laarhoven, P.J.M. (1985), Parallel variable metric algorithm for unconstrained optimization, {\it Mathematical Programming}, {\bf 33}, 68--87.


\end{references}
\end{document}